\documentclass[12pt,a4paper]{article}
\usepackage{amssymb}

\usepackage{graphicx,amssymb,amsfonts,epsfig,amsthm,a4,amsmath,url}
\usepackage[latin1]{inputenc}

\newtheorem{ThmIntro}{Theorem}
\newtheorem{CorIntro}[ThmIntro]{Corollary}

\newtheorem{thm}{Theorem}[section]
\newtheorem{cor}[thm]{Corollary}
\newtheorem{lem}[thm]{Lemma}

\newtheorem{prop}[thm]{Proposition}
\theoremstyle{definition}
\newtheorem{defn}[thm]{Definition}

\theoremstyle{remark}
\newtheorem{rem}[thm]{Remark}
\theoremstyle{cremark}
\newtheorem{crem}[thm]{Crucial remark}

\newtheorem{ex}[thm]{Example}

\numberwithin{equation}{section}

\newcommand{\N}{\mathbf{N}}
\newcommand{\R}{\mathbf{R}}

\newcommand{\Isom}{\text{Isom}}

\newcommand{\bpr}{\noindent \textbf{Proof}: ~}

\newcommand{\epr}{~$\blacksquare$}

\newcommand{\eps}{\varepsilon}

\newcommand{\Aut}{\textnormal{Aut}}

\newcommand{\A}{\mathcal{A}}

\title{Large scale Sobolev inequalities on metric measure spaces and applications.}
\author{Romain Tessera}
\date{\today}

\begin{document}

\baselineskip=16pt

\maketitle

\begin{abstract}
We introduce a notion of ``gradient at a given scale" of functions
defined on a metric measure space. We then use it to define Sobolev
inequalities at large scale and we prove their invariance under
large-scale equivalence (maps that generalize the quasi-isometries).
We prove that for a Riemmanian manifold satisfying a local
Poincar\'e inequality, our notion of Sobolev inequalities at large
scale is equivalent to its classical version. These notions provide
a natural and efficient point of view to study the relations between
the large time on-diagonal behavior of random walks and the
isoperimetry of the space. Specializing our main result to locally
compact groups, we obtain that the $L^p$-isoperimetric profile, for
every $1\leq p\leq \infty$ is invariant under quasi-isometry between
amenable unimodular compactly generated locally compact groups. A
qualitative application of this new approach is a very general
characterization of the existence of a spectral gap on a
quasi-transitive measure space $X$, providing a natural point of
view to understand this phenomenon.
\medskip

    \hfill\break \noindent {\sl Mathematics Subject
Classification:} Primary 51F99; Secondary 43A85. \hfill\break {\sl
Key words and Phrases:} Sobolev inequalities, isoperimetry,
analysis on metric measure spaces, coarse geometry, random walks,
spectral gap, locally compact groups.

\end{abstract}

\tableofcontents

\section{Introduction}
We introduce a notion of ``gradient at a certain scale" of a
bounded function defined on a general metric measure space. We
then give a meaning to the notion of ``large-scale" Sobolev
inequalities\footnote{This functional analysis approach
generalizes the purely geometric notion of large-scale
isoperimetry that we introduced in \cite{tess'}.} for metric
measure spaces and we show their invariance under large-scale
equivalence. Moreover, we show that under some controlled
connectivity assumption, the large scale Sobolev inequalities are
equivalent to Sobolev and inequalities at a positive given scale.
We also study the relations between our notion of gradient at
given scale and the well-known infinitesimal notion of generalized
upper-gradient. In particular, we prove that for a riemannian
manifold satisfying a local Poincaré inequality, our large-scale
Sobolev inequalities are equivalent to their usual versions
(defined with the riemannian gradient). The improvement of our
point of view is to get rid of any condition at small scale since
it is rubbed out by the definition of the large-scale gradient.
This level of generality can be really useful, for instance for
the study of $\sigma$-compact locally compact groups where no nice
local structure is available. It can be also important to include
(highly) non-geodesic spaces, as subspaces of a metric space are
not coarsely geodesic in general (this can be the case of a
subgroup equipped with the induced distance). Moreover, note that
a locally compact group has no coarsely geodesic left invariant
proper metric unless it is compactly generated (see
Proposition~\ref{connectivitegoupprop}).

These notions provide a natural and efficient point of view to study
the relations between the large time on-diagonal behavior of random
walks and the isoperimetry of the space. In particular, we obtain
that, under mild assumptions on a metric measure space, upper bounds
on the probability of return of symmetric random walks are
characterized by large-scale Sobolev inequalities, and therefore are
invariant under large-scale equivalence (see Theorem~\ref{randomthm}
for a precise statement).

As a qualitative application, we prove that a reversible random walk
on a quasi-transitive measure space has spectral radius equal to $1$
if and only if the group acting is amenable and unimodular. This
provides a general and direct explanation for a phenomenon that has
been proved in particular cases\footnote{Note that some of the
results of these articles are more precise than ours and in a sense,
more general when they manage to deal with non-reversible random
walk.} in \cite{Kest,Bro,Salv,SoW,P',SaW}.

\medskip

\subsection*{Statement of the main results in the homogeneous setting}

Let us present present our results in a very special --though
interesting-- case: when $X=G$ is a group. Let $G$ be a locally
compact, compactly generated group equipped with a left-invariant
Haar measure $\mu$. Let $S$ be a compact subset of $G$ such that
$\bigcup_{n\in \N}S^n=G$. Equip $G$ with the left-invariant word
metric\footnote{To obtain a real metric, one must assume that $S$
is symmetric, but this does not really play a role in the sequel.}
associated to $S$, $d_S(g,h)=\inf\{n,\; g^{-1}h\in S^n\}$.

\subsubsection*{Quantitative results}

Recall that a quasi-isometry between two metric spaces $(X,d_X)$ and
$(Y,d_Y)$ is a map $F:\:X\to Y$ which is bi-Lipschitz for large
distances, i.e.
$$C^{-1}d_X(x,y)-C\leq d_Y(F(x),F(y))\leq Cd_X(x,y)+C,$$
for any $x,y\in X$, $C$ being a positive constant; and almost
surjective, i.e.
$$\sup_{z\in Y} d(z,F(X))<\infty.$$

\medskip

Let $\lambda$ be the action of $G$ by left-translations on functions
on $G$, i.e. $\lambda(g)f(x)=f(g^{-1}x)$. For any $1\leq p\leq
\infty$, and any subset $A$ of $G$, define
$$J_p(A)=\sup_{f}\frac{\|f\|_p}{\sup_{s\in S}\|f-\lambda(s)f\|_p},$$
where $f$ runs over elements of the space $L^p(A)$
($L^p$-functions supported in $A$).

We can define two kinds of ``$L^p$-isoperimetric profile",
depending on whether we want to optimize $J_p(A)$ fixing the
volume of $A$, or its diameter. In the first case, we obtain what
is often called the $L^p$-isoperimetric profile (see for instance
\cite{C,Coulhon}),
$$j_{G,p}(v)=\sup_{\mu(A)=v}J_p(A).$$
In the second case, we obtain what we call the $L^p$-isoperimetric
profile inside balls since it is given by
$$J^b_{G,p}(n)=\sup_{x\in X} J_p(S^n).$$

\medskip

We will be interested in the ``asymptotic behavior" of these
nondecreasing functions. Precisely, let $f,g: \R_+\to\R_+$ be
nondecreasing functions. We write respectively $f\preceq g$, $f\prec
g$ if there exists $C>0$ such that $f(t)=O(g(Ct))$, resp.
$f(t)=o(g(Ct))$ when $t\to \infty$. We write $f\approx g$ if both
$f\preceq g$ and $g\preceq f$. The asymptotic behavior of $f$ is its
class modulo the equivalence relation $\approx$.

\medskip

Now, we can state our main results in this setting.
\begin{ThmIntro}[\textnormal{see Corollary~\ref{groupisocor1}}]
Assume that $(G,S)$ and $(H,T)$ are two unimodular compactly
generated, locally compact groups, equipped with symmetric
generating subsets $S$ and $T$ respectively. Then, the asymptotic
behaviours of $j_{G,p}$, $J^b_{G,p}$, for any $1\leq p\leq \infty$
do not depend on $S$. Moreover, if $G$ is quasi-isometric to $H$,
then
$$j_{G,p}\approx j_{H,p},$$ and
$$J^b_{G,p}\approx J^b_{H,p}.$$
\end{ThmIntro}

\subsubsection*{A qualitative result}

We also derive a qualitative result on quasi-transitive spaces. Let
$G$ be a locally compact, compactly generated group. Let $(X,\mu)$
be a quasi-transitive $G$-space, i.e. a locally compact Borel
measure space on which $G$ acts measurably, co-compactly, properly,
and almost preserving the measure $\mu$, i.e.
$$\sup_{g\in G}\sup_{x\in X}\frac{d(g\cdot\mu)}{d\mu}(x)<\infty.$$
For every $x\in X$, let $\nu_x$ be a probability measure on $X$
which is absolutely continuous with respect to $\mu$. We assume
that there exist $S\subset S'$, two compact generating subsets of
$G$, and a compact subset of $X$ satisfying $GK=X$, such that for
every $x\in X$, the support of $\nu_x$ is contained in $gS'K,$ for
some $g\in G$ such that $x\in gSK$. Let us also suppose that
$\nu_x(y)$ is larger than a constant $c>0$ for $y$ in $gSK$.
Denote by $P$ the Markov operator on $L^2(X)$ defined by
$$Pf(x)=\int f(gy)d\nu_x(y).$$
We make the (important) assumption that $P$ is self-adjoint.
\begin{ThmIntro}\label{mainThmIntro}\textnormal{(see Theorem \ref{quasitransitivethm} and Corollary~\ref{nongeomCheegerconstant})}
The following are equivalent
\begin{itemize}\item the spectral radius of $P$ is less than $1$;
\item $G$ is either non-unimodular or non-amenable.\item $G$ is quasi-isometric to
a graph (of bounded degree) with positive Cheeger constant.
\end{itemize}
\end{ThmIntro}
This theorem is a slight generalization of the following recent
result of Saloff-Coste and Woess \cite{SaW}, which they obtained by
completely different (and less elementary) methods.

\begin{CorIntro}\textnormal{\cite{SaW}}
Let $(X,d)$ be a geodesic metric space and let $G$ be a closed
subgroup of $\Isom(M)$ acting co-compactly on $X$ by isometries. Fix
$r>0$. Then $G$ is unimodular and amenable if and only if the
spectral radius of the average operator on balls of radius $r$ is
$1$.
\end{CorIntro}

Our approach unifies the following results, enhancing their
``large-scale" nature. An obvious particular case is when the space
$X$ is the group itself.

\begin{CorIntro}
Let $G$ be a locally compact group equipped with a left Haar measure
$\mu$. Then $G$ is unimodular and amenable if and only if for every
compactly supported, symmetric (with respect to $\mu)$ random walk
on $G$ has spectral radius $1$.
\end{CorIntro}

\begin{CorIntro}\textnormal{\cite{Salv}}
Let $X$ be a connected graph of bounded degree and let $G$ be a
closed subgroup of $\Aut(X)$ such that $X/\Aut(X)$ is finite. Then
$G$ is unimodular and amenable if and only if the spectral radius of
the simple random walk equals $1$.
\end{CorIntro}
When $G$ is transitive this theorem has been proved in \cite{SoW}.

\begin{CorIntro}\textnormal{\cite{SaW}}\textnormal{(see Corollary~\ref{riamannianSpectralGap})}
Let $M$ be a Riemannian\footnote{Actually the authors give a method
that allow them to treat a large class of examples, like all the
examples given here, excepted the case of the group itself as they
need $X$ to be geodesic.} manifold and let $G$ be a closed subgroup
of $\Isom(M)$ acting co-compactly on $M$. Then $G$ is unimodular and
amenable if and only if the spectral radius of the heat kernel on
$M$ equals $1$, or in other words, if the (Riemannian) Laplacian on
$M$ has no spectral gap around zero.
\end{CorIntro}
The case where $G$ is transitive has been treated in \cite{P'} and
the case where $M$ is the universal cover of a compact manifold has
been proved in \cite{Bro}.

\medskip

\subsection*{Organization of the paper}
\begin{itemize}
\item In Section~\ref{gradientDefSection}, we introduce a notion
of Sobolev inequalities that capture the geometry at a scale larger
than $h>0$.

\item In Section~\ref{SoboInterpretationSection}, we discuss the
geometric and probabilistic interpretations of those Sobolev
inequalities.

\item In Section~\ref{Sobo/profileSection}, we discuss the relations
between Sobolev inequalities and the isoperimetric profile.

\item In Section~\ref{largescaleDefSection}, we introduce the notion
of large-scale equivalence, which is a measured version of the
well-known notion of coarse equivalence (see \cite{Ro}).

\item In section~\ref{Examples}, we discuss some examples of
large-scale equivalences in the contexts of locally compact groups,
manifolds, graphs etc.

\item In Section~\ref{equivalentgradientsection}, we prove a
technical but important fact: the definition of large-scale Sobolev
inequalities does not depend on the choice of a ``large-scale"
gradient.

\item In Section~\ref{largescaleinvariancesection}, we prove our
main result, namely that large-scale Sobolev inequalities are
invariant under large-scale equivalence.

\item In Section~\ref{random}, we relate Sobolev inequalities to
upper bounds on the probability of return of symmetric random
walks.

\item In Section~\ref{controlscalesection}, we discuss the validity
of Sobolev inequalities at a given scale when it is true at
large-scale. In particular, in
Sections~\ref{finite->infinitesimalsection} and
\ref{infinitesimal->finitesection}, we prove that under some mild
local assumptions, the large-scale Sobolev inequalities are
equivalent to their classical versions on a Riemaniann manifold.

\item Finally, in Section~\ref{HomogeneousSection}, we prove the
results announced in the introduction in the context of locally
compact groups and quasi-transitive spaces.
\end{itemize}

\section{Functional analysis at a given scale}\label{gradientDefSection}

\subsection{Local norm of gradient at scale
$h$}\label{largescalegradient} Let $(X,d)$ be a metric space. The
purpose of this section is to define a notion of gradient that
capture the geometry at a certain scale --say $h$-- of $X$. More
precisely, as we will see in sequel, what we really need to define
is not the gradient of a function itself, but rather a local norm of
this gradient (that plays the role of the modulus of the gradient
for a Riemannian manifold). The first naive idea to do this is to
define
$$|\nabla f|_h(x)=\sup_{y\in B(x,h)}|f(y)-f(x)|$$
for any function $f\in L^{\infty}(X)$, $B(x,h)$ denoting the closed
ball of center $x$ and radius $h$. Note that this can be written in
the following form:
$$|\nabla f|_h(x)=\|f-f(x)\|_{\infty,B(x,h)}$$
which emphasizes the fact that we actually consider a ``local"
$L^{\infty}$-norm. Naturally, we would like to define also the local
$L^p$-norm of the gradient of $f$, for every $1\leq p\leq \infty$.
For this, we obviously need a measure on $X$. What we could do is
start from a measure on $X$ and define a local $L^p$-norm as the
$L^p$ norm restricted to balls with respect to this measure.
However, when we consider a random process on $X$, the notion of
local $L^2$-norm that naturally emerges is the $L^2$-norm with
respect to the probability transition. This motivates the following
definition.




Let $(X,d,\mu)$ be a metric measure space. Consider a family
$P=(P_x)_{x\in X}$ of probability measures on $X$. Then for every
$p\in [1,\infty]$, we define an operator $|\nabla|_{P,p}$ on $
L^{\infty}(X)$ by
$$\forall f\in  L^{\infty}(X),\quad |\nabla f|_{P,p}(x)=\| f-f(x)\|_{P_x,p}=\left(\int |f(y)-f(x)|^pdP_x(y)\right)^{1/p},$$
if $p<\infty$; and for $p=\infty$, we decide that
$$|\nabla f|_{P,\infty}(x)= \|
f-f(x)\|_{P_x,\infty}=\sup\{|f(y)-f(x)|, \; y\in Supp(P_x)\}.$$

\begin{defn} A family of probabilities $P=(P_x)_{x\in X}$ on $X$ is called a
viewpoint at scale $h>0$ on $X$ if there exist a large constant
$1\leq A<\infty$ and a small constant $c>0$ such that for
($\mu$-almost) every $x\in X$:
\begin{itemize}
\item $P_x\ll\mu;$

\item $p_x=dP_x/d\mu$ is supported in $B(x,Ah);$

\item $p_x$ is larger than $c$ on $B(x,h)$.
\end{itemize}
\end{defn}

\begin{rem}
Note that a viewpoint at scale $h$ is also a viewpoint at scale $h'$
for any $h'<h.$
\end{rem}

\begin{ex}
A basic example of viewpoint at scale $h$ is given by
$$P_x=\frac{1}{V(x,h)}1_{B(x,h)},\quad \forall x\in X,$$
where $V(x,r)$ denotes the volume of the ball centered at $x$ of
radius $r$. We denote the associated $L^p$-gradient by
$|\nabla|_{h,p}$. Note that with the notation of the beginning of
the Section~,
$$|\nabla|_h=|\nabla|_{h,\infty}.$$
\end{ex}

\begin{ex}\label{graphEx}(For more details, see \cite{C})
To any connected simplicial graph, we associate a metric measure
space, whose elements are the vertices of the graph, the measure is
the counting measure and the distance is the usual discrete geodesic
distance for which two distinct points joined by an edge are at
distance $1$ from one another. For simplicity, we will simply call
such a metric measure space a graph. The usual discrete local norm
of gradient on a graph, usually denoted by $|\nabla f (x)|$,
corresponds\footnote{In \cite{C}, they consider a slightly different
definition, where the average is taken over the set of neighbors of
$x$ instead of the ball $B(x,1)$.} with our notations to
$$|\nabla f|_{1,2}(x)=\left(\frac{1}{V(x,1)}\sum_{y\in B(x,1)}|f(x)-f(y)|^2\right)^{1/2}.$$
\end{ex}

\begin{rem}{\bf (Interpretations of the notion of viewpoint at scale $h$.) }
A viewpoint at scale $h$ has at least two interesting
interpretations: one as an operator transition of a random walk on
$X$; the other as a Markov operator acting on $L^p(X)$ for every
$p\geq 1$. This operator is defined by
$$Pf(x)=\int_{X}f(y)dP_x(y).$$
Consequently, there is a natural semi-group structure on the set of
viewpoints at scale $h$ on space $X$. Indeed, it is straightforward
to check\footnote{One has to suppose that the space is locally
doubling: see Definition~\ref{doublinganyscaledefn}.} that if $P$ is
a viewpoint at scale $h$ and $Q$ is a viewpoint at scale $h'$, then
$P\circ Q$ is a viewpoint at any scale $h"<h+h'$.
\end{rem}

\begin{rem}\label{gradientveritablerem}{\bf (Alternative definition of gradient at scale $h$.) }
Let us indicate another way of describing the objects that we
introduced. Instead of directly defining a {\it local norm  of the
gradient at scale $h$}, we could first define a true {\it gradient
at scale $h$} on a fiber space over $X$ and then take a local norm
of the gradient on the fibers. Here the fiber space would be
$Y_h=\{(x,y)\in X^2, d(x,y)\leq h\}$ with projection $\pi: Y\to X$
on the first factor, so that $\pi^{-1}(x)$ identifies with
$B(x,h)$. The gradient at scale $h$ of $f$ is then $\nabla_h
f(x,y)=f(x)-f(y)$, where $(x,y)\in Y_h$. A viewpoint at scale $h$
on $X$ is now a probability measure on every fiber of some
$Y_{Ah}$ for $A$ large enough; and the $L^p$-gradient of $f$
associated to such a viewpoint corresponds to the $L^p$-norm of
$f$ in every fiber with respect to this measure\footnote{Note that
we can also define the gradient of $f$ without referring to the
scale: $\nabla f: X\times X\to \R,$ $\nabla f(x,y)=f(x)-f(y)$,
looking at $X\times X$ as a fiber space over the first factor.
Then the scale appears when choosing a norm on every fiber.}.
\end{rem}

\subsection{Laplacian at
scale $h$}\label{laplacianSection}

We can also define a Laplacian w.r.t. a viewpoint $P=(P_x)_{x\in X}$
by
$$\Delta_{P}f(x)=(id-P)f(x),$$ and more generally a $p$-Laplacian
for any $p>1$ by
$$\Delta_{P,p}f(x)=\int |f-f(x)|^{p-2}(f-f(x))dP_x.$$
If $P$ is self-adjoint with respect to the scalar product associated
to $\mu$, then we have the usual relations
$$\langle \Delta_{P,p}f,g\rangle=\int \left(\int |f(y)-f(x)|^{p-2}(f(y)-f(x))(g-g(x))dP_x(y)\right)d\mu(x),$$
$$\langle \Delta_{P,p}f,f\rangle=\int\int |f(y)-f(x)|^{p}p_x(y)d\mu(y)d\mu(x)=\||\nabla f|_{P,p}\|_p^p,$$
and in particular, for $p=2$,
$$\langle \Delta_{P}f,f\rangle=\int\int |f(y)-f(x)|^2p_x(y)d\mu(x)d\mu(y)=\||\nabla f|_{P,2}\|_2^2.$$

In particular, if $A$ is a measurable subset of $X$, The first
eigenvalue $\delta_P$ of $\Delta_{P}$ acting on square-integrable
functions supported by $A$ is $$\delta_P(A)=\inf_f
\frac{\||\nabla_{P,2} f|\|_2^2}{\|f\|_2^2},$$ where $f$ runs over
square-integrable functions supported by $A$.

\subsection{Sobolev inequalities at scale $h$}
 Let $\varphi: \R_+\to\R_+$ be an increasing function and let
$p\in [1,\infty]$. The following formulation of Sobolev inequality
was introduced in \cite{Cou}. We refer to \cite{Coulhon} for the
link with more classical formulations, for instance in $\R^n$.

\begin{defn}\label{sobolevdef} One says that $X$ satisfies a Sobolev inequality
$(S_{\varphi}^p)$ at scale (at least) $h>0$ if there exists some
finite positive constants $C$, $C'$ depending only on $h$, $p$ and
$\varphi$ such that
$$\| f\|_p\leq C\varphi(C'|\Omega|)\| |\nabla f|_h \|_p$$
where $\Omega$ ranges over all compact subsets of $X$, $|\Omega|$
denotes the measure $\mu(\Omega)$, and $f\in L^{\infty}(\Omega)$, $
L^{\infty}(\Omega)$ being the set of elements of $L^{\infty}(X)$
with support in $\Omega$.
\end{defn}


\begin{defn}
We say that $X$ satisfies a large-scale Sobolev inequality
$(S_{\varphi}^p)$ if it satisfies $(S_{\varphi}^p)$ at some scale
$h$ (equivalently, for $h$ large enough).
\end{defn}
\begin{crem}
Note that to define the Sobolev inequalities at large scale, we
arbitrarily chose to write them with $|\nabla|_h$ whereas we could
have defined them with $|\nabla|_{P,q}$ for any viewpoint
$(P_x)_{x\in X}$ at scale $h$ and any $q\geq 1$. A crucial and
useful fact that we prove in Section~\ref{equivalentgradientsection}
is that satisfying a {\it large-scale} Sobolev inequality {\bf does
not} depend on this choice.
\end{crem}

\begin{rem}
Note that for large scale Sobolev inequalities, only $\Omega$ with
large volume are involved. In fact, we will only be interested in
the asymptotic behavior of $\varphi$.
\end{rem}

\begin{rem}\label{Sobopq}
It is easy to prove that $(S_{\varphi}^p)$ implies $(S_{\varphi}^q)$
whenever $p\leq q<\infty$ for any choice of gradient (see
\cite{Coulhon} for a proof in the Riemannian setting). It is proved
in \cite{CL} that the converse is false for general Riemannian
manifolds. This is likely to be true for groups, although it is
still open.
\end{rem}

\subsection{Link with Sobolev inequalities for infinitesimal
gradients}\label{soboscalesection}

Other notions of ``local norm of gradient" have been introduced and
studied for general metric spaces. In particular the notion of upper
gradient plays a crucial role in the study of doubling metric spaces
equipped with the Hausdorff measure (see for instance
\cite{Hei,Sem}, or Definition~\ref{uppergradDef}). Those spaces
naturally occur as boundaries of Gromov-hyperbolic spaces and are
often studied up to quasi-conformal maps. Such a point of view is
quite different from ours since it focuses on the local properties
of the space, which is often supposed compact. However, it is
natural to ask when a Sobolev inequality at large scale is
equivalent to the same Sobolev inequality w.r.t. some upper
gradient. In particular, given a Riemannian manifold, is it true
that it satisfies a Sobolev inequality at large scale if and only if
it satisfies it for its usual gradient?
Proposition~\ref{infinitesimal->fini} says that if a Riemannian
manifold satisfies a Sobolev inequality for its usual gradient, then
it also satisfies it at large scale (but the proof is not as obvious
as one could expect). However, the converse can be false, for
instance if the Riemannian manifold contains a sequence of open
submanifolds isometric to open half-spheres of radius going to zero.
A sufficient condition to get a positive answer is to ask for a
local Poincaré inequality (see
Proposition~\ref{finite->infinitesimal}).

\medskip

\noindent{\bf Other ideas for ignoring the local geometry of a
Riemannian manifold.} Different strategies have been used to ignore
the local geometric properties of a manifold. In \cite{ChaFel} for
instance, they avoid the local behavior of the isoperimetric profile
on a manifold by restricting it to subsets containing a geodesic
ball of fixed radius. In \cite{Coulhon'}, they consider Nash
inequalities restricted to functions convoluted by the heat kernel
at time $1$ and obtain in this way the invariance under
quasi-isometries of certain upper bounds of the on-diagonal
behaviour of the heat kernel: this idea is quite closed to ours (see
Remark~\ref{finalrem}). This issues are discussed in
Sections~\ref{infinitesimal->finitesection} and
\ref{finite->infinitesimalsection}. Among other things, we prove
under a very weak property of bounded geometry that a manifold
satisfies a Sobolev inequality at large scale if and only if it
satisfies it for the usual gradient in restriction to functions of
the form $g=Pf$, where $P$ is the Markov operator associated to any
viewpoint at some scale $h>0$.

\section{Sobolev inequalities $(S_{\varphi}^p)$ at scale $h$ for
$p=1,2,\infty$}\label{SoboInterpretationSection}

Now let us give characterizations of $(S_{\varphi}^p)$ at given
scales for the important values of $p=1,2,\infty$ (see \cite{C} for
the case of graphs and \cite{Coulhon} for Riemannian manifolds).

\subsection{$(S_{\varphi}^{\infty})$ and volume growth}

In \cite{Coulhon'} (see also \cite[proposition~22]{Coulhon}), it
is proved that $(S_{\varphi}^{\infty})$ can only hold if $\varphi$
is unbounded and then is equivalent to the volume lower bound
$$V(x,r)\geq \varphi^{-1}(r)$$
where $\varphi^{-1}(r)=\{v,\varphi(v)\geq r\}$, for every $x\in X$
and every $r>0$. The original proof works formally in our setting.

\begin{prop}
Let $(X,d\mu)$ be a metric measure space. The Sobolev inequality
$(S_{\varphi}^{\infty})$ at scale $h$ can only hold if $\varphi$ is
unbounded and then is equivalent to the volume lower bound
$$V(x,r)\geq \varphi^{-1}(r)$$
for $r\geq h$.~$\blacksquare$
\end{prop}

\subsection{$(S_{\varphi}^1)$ and isoperimetry} \label{geominterpretation}
\begin{prop}\label{isoperimetryProp}
The inequality $(S_{\varphi}^1)$ at scale $h$ is equivalent to the
isoperimetric inequality (at scale $h$)
$$\frac{|\partial_h \Omega|}{|\Omega|}\geq \frac{1}{C\varphi(C'|\Omega|)}$$
where the boundary of $A$ is defined by $$\partial_h A=[A]_h\cap
[A^c]_h$$ with the usual notation $[A]_h=\{x\in X,\; d(x,A)\leq
h\}.$\epr
\end{prop}
The usual proof of this equivalence (see \cite{Coulhon}) works
formally in our context, using the following version of the
co-area formula.
\begin{lem}{\bf (Co-area formula at scale $h$)}
\begin{equation}\label{coaire}
\frac{1}{2}\int_{\R_+}\mu\left(\partial_h \{f\geq t\}\right)dt\leq
\int_{X}|\nabla f|_h(x)d\mu(x)\leq\int_{\R_+}\mu\left(\partial_h
\{f\geq t\}\right)dt
\end{equation}
where $f$ is a non-negative measurable function defined on $X$.
\end{lem}
\bpr For every measurable subset $A\subset X$, we have
$$\mu(\partial_h A)=\int_X|\nabla 1_A|_h(x)d\mu(x).$$
Thus, (\ref{coaire}) follows by integrating over $X$ the following
local inequalities
\begin{equation}\label{localcoaire}
\frac{1}{2}\int_{\R_+}|\nabla 1_{\{f\geq t\}}|_h(x) dt\leq |\nabla
f|_h(x)\leq\int_{\R_+}|\nabla 1_{\{f\geq t\}}|_h(x) dt,
\end{equation}
for every $x\in X.$ The right-hand inequality results from the fact
that $f=\int_{\R_+}1_{\{f\geq t\}}dt$ and from the sub-additivity of
$|\nabla|_h$. To prove the left-hand, note that $|\nabla 1_{\{f\geq
t\}}(x)|_h=1$ if and only if
$$\inf_{B(x,h)}f <t \leq \sup_{B(x,h)}f$$
or
$$\inf_{B(x,h)}f\leq t < \sup_{B(x,h)}f;$$
Hence, $$\int_{\R_+}|\nabla 1_{\{f\geq t\}}|_h(x) dt\leq
\sup_{B(x,h)}f-\inf_{B(x,h)}f\leq 2|\nabla f|_h(x),$$ which proves
(\ref{localcoaire}).\epr

\subsection{Probabilistic interpretation of $(S_{\varphi}^2)$ at
scale $h$}

The case $p=2$ is of particular interest since it contains some
probabilistic information on the space $X$. It is proved in
\cite{CG} that for manifolds with bounded geometry, upper bounds of
the large-time on-diagonal behavior of the heat kernel are
equivalent to some Sobolev inequality $(S_{\varphi}^2)$. In
\cite{C}, a similar statement is proved for the standard random walk
on a weighted graph. In Section~~\ref{random}, we give a
discrete-time version of this theorem in our general setting. The
proof of Theorem~\ref{randomthm} below emphasizes the fact that the
notion of viewpoint at scale $h$ is likely to be the most natural
way of capturing the link between large-scale geometry and the
long-time behavior of random walks on $X$.
\begin{defn}
Let $(X,d,\mu)$ be a metric measure space and consider some $h>0$. A
view-point $P=(P_x)_{x\in X}$ at scale $h$ on $X$ is called
symmetric if one of the following equivalent statement holds.
\begin{itemize}
\item The random walk whose probability of transition is $P$ is
reversible with respect to the measure $\mu$.

\item The associated operator on $L^2(X,\mu)$ defined by
$$Pf(x)=\int_X f(y)dP_x(y)$$
is self-adjoint.

\item For every a.e. $x,y\in X$, $p_x(y)=p_y(x).$

\end{itemize}
\end{defn}

\begin{defn}
We call a reversible random walk at scale $h$ a random walk whose
probability transition is a symmetric view-point at scale $h$.
\end{defn}

\begin{ex}
Let $(X,d,\mu)$ be a metric measure space. Consider the standard
viewpoint at scale $h$ of density $p_x=1_{B(x,h)}/V(x,h)$ with
respect to $\mu$. In general, this is not a symmetric viewpoint,
i.e. the random walk of probability transition
$dP_x(y)=p_x(y)d\mu(y)$ is not reversible with respect to $\mu$.
However, it is reversible with respect to the measure $\mu'$ defined
by
$$d\mu'(x)= V(x,h)d\mu(x).$$
It is easy to check that if $(X,d,\mu)$ is locally doubling, then so
is $(X,d,\mu')$. Moreover, if $x\mapsto V(x,h)$ is bounded from
above and from below, then $P$ defines a symmetric viewpoint on
$(X,d,\mu')$.
\end{ex}

The relations between large-scale Sobolev inequalities
$(S_{\varphi}^2)$ and random walks on a metric measure space are
summarized in the following theorem, whose proof is adapted from
\cite[Theorem~7.2]{Coulhon}. We use the notation
$dP_x^n(y)=p_x^{n}(y)d\mu(y)$.

\begin{thm}\label{randomthm}\textnormal{(see Section~\ref{random})}
Let $X=(X,d,\mu)$ be a metric measure space and let $P=(P_x)_{x\in
X}$ be a symmetric view-point at scale $h$ on $X$. Let $\varphi$ be
some increasing positive function. Define $\gamma$ by
$$t=\int_{0}^{1/\gamma(t)}\left(\varphi(v)\right)^2\frac{dv}{v}.$$
\begin{itemize}
\item[(i)] Assume that $X$ satisfies a large-scale Sobolev inequality
$(S_{\varphi}^2)$. Then
$$p_x^{2n}(x)\leq \gamma(cn) \quad \forall n\in \N, a.e \forall x\in X,$$
for some constant $c>0$.

\item[(ii)] If the logarithmic derivative of $\gamma$ has at most polynomial
growth\footnote{This condition, called $(\delta)$ in \cite[p~18]{C}
is very weak since it is satisfied by all functions $(\log
t)^{a}t^{b}e^{ct^d}$ for any real numbers $a,b,c,d$.} and if
$$p_x^{2n}(x)\leq \gamma(n) \quad \forall n\in \N, a.e \forall x\in X,$$
then $X$ satisfies $(S_{\varphi}^2)$ w.r.t. $|\nabla|_{P,2}$.
\end{itemize}
\end{thm}


\section{Sobolev inequalities and $L^p$-isoperimetric profiles at scale
$h$}\label{Sobo/profileSection}

\subsection{$L^p$-isoperimetric profile at scale
$h$}\label{isoperimetricprofilesection}

Generalizing the case $p=1$, Sobolev inequalities $(S_{\varphi}^p)$
can be also understood as $L^p$-isoperimetric inequalities. Let $A$
be a measurable subset of $X$. For every $p>0$, define
$$J_{p}(A)=\sup_{f} \frac{\| f\|_p}{\||\nabla f|_h \|_p}$$
where the supremum is taken over functions $f\in L^{\infty}(A)$.
Note that for $p=2$, this is just the square root of the inverse of
the first eigenvalue of the Laplacian $\Delta_P$ acting on
square-integrable functions supported by $A$ (see
Section~\ref{laplacianSection}).

Now, taking the supremum over subsets $A$ with measure less than
$m>0$, we get an increasing function $j_{X,p}$ sometimes called
the $L^p$-isoperimetric profile. Clearly, the space $X$ always
satisfies the Sobolev inequality $(S_{\varphi}^p)$ with $\varphi=
j_{X,p}$. Conversely, if $X$ satisfies $(S_{\varphi}^p)$ for a
function $\varphi$, then
$$j_{X,p}\succeq \varphi.$$ It is easy to check that
$$j_{X,p}\preceq j_{X,q}$$
whenever $p\leq q<\infty$ (see Remark~\ref{Sobopq} about Sobolev
inequalities).

Note that the terminology ``isoperimetric profile" is somewhat
ambiguous since there exist various nonequivalent definitions (see
in particular \cite[Chapter 1]{CS}). One of them is
$$j_X(m)=\sup_{|A|\leq m}\frac{|A|}{|\partial_h A|}.$$
As a corollary of Proposition~\ref{isoperimetryProp}, we have
\begin{prop}
The $L^1$-isoperimetric profile and the isoperimetric profile have
the same asymptotic behavior, i.e. $$j_X\approx j_{X,1},$$ taking
the same $h$ in the definition of the gradient and in the
definition of the boundary.
\end{prop}

\subsection{$L^p$-isoperimetric profile inside balls}

\begin{defn}\label{isoperimetricprofiledef}
Let us fix a gradient at scale $h$ on $X$. The $L^p$-isoperimetric
profile inside balls is the nondecreasing function $J^b_{G,p}$
defined by
$$J^b_{X,p}(t)=\sup_{x\in X}J_p(B(x,t)).$$
\end{defn}
Note that $J^b_{X,p}(t)$ is the supremum of $J_p(A)$ over subsets
$A$ of diameter\footnote{This profile is associated to another kind
of Sobolev inequalities, where the function $\varphi$ of the volume
is replaced by a function $\Phi$ of the diameter.} less than $t$.
The $L^p$-isoperimetric profile inside balls plays a crucial role in
the study of uniform embeddings of amenable groups into $L^p$-spaces
(see~\cite{tessera}). It is also central in the proof
\cite{tessera'} that a closed at infinity, homogenous manifold does
not carry any non-constant $p$-harmonic function with gradient in
$L^p$.

\subsection{Link with the large-scale isoperimetry introduced
in \cite{tess'}}

One can also define another kind of isoperimetric profile at scale
$h$:
$$I(t)=\inf_{\mu(A)\geq t}\mu(\partial_h A)$$
which can be specialized on a family of (measurable) subsets of
finite volume $\A$: we call lower (resp. upper) profile at scale $h$
restricted to $\A$ the nondecreasing function $I^{\downarrow}_{\A}$
defined by
$$I^{\downarrow}_{\A}(t)=\inf_{\mu(A)\geq t, A\in \A}\mu(\partial_h A)$$
(resp. $I^{\uparrow}_{\A}(t)=\sup_{\mu(A)\leq t, A\in
\A}\mu(\partial_h A)$). We can then study the large scale
isoperimetric properties of a family $\A$ considering the asymptotic
behavior of these two increasing functions \cite{tess'}. In
\cite{tess'}, we used this variant to investigate the question: are
balls always asymptotically isoperimetric in a metric measure space
with doubling property? For that purpose, we introduced a general
setting adapted to the study of asymptotic isoperimetry on metric
measure spaces. An important consequence of the geometric
interpretation of Sobolev inequalities in $L^1$ (see
Section~~\ref{geominterpretation}) is that every geometric notion
that we introduced in \cite[Section~3]{tess'} appears as a
particular case of the functional point of view adopted in the
present paper. In particular, \cite[Theorem~3.10]{tess'} that
implied the invariance under large-scale equivalence of
isoperimetric properties is now covered by the lemmas of
Section~~\ref{largescaleinvariancesection}. Moreover, we choose here
to treat separately the large-scale setting, where no connectivity
hypotheses are required on the spaces, and the control on the scale
that really {\it depends} on a connectivity assumption (see
Section~~\ref{controlscalesection}).


\section{Large-scale equivalence between metric measure
spaces}\label{largescaleDefSection}

In this section, we define an equivalence relation, called
large-scale equivalence between metric measure spaces. This notion
is simply an adaptation of the notion of coarse equivalence for
metric spaces introduced by Roe in \cite{Ro}, for spaces endowed
with a measure.

The metric measure spaces that we will consider satisfy a very weak
property of bounded geometry introduced in \cite{CS}.
\begin{defn}\label{doublinganyscaledefn}
We say\footnote{In \cite{CS} and in \cite{tess'}, the  local
doubling property is denoted $(DV)_{loc}$.} that a space $X$ is
locally doubling at scale $r>0$ if there exists a constant $C_r$
such that
$$\forall x\in X, \quad V(x,2r)\leq C_rV(x,r)$$ where $V(x,r)=\mu(B(x,r)).$
If it is locally doubling at every scale $r>0$, then we just say
that $X$ is locally doubling.
\end{defn}
\begin{rem}
Since the constant $C_r$ depends on $r$, the locally doubling
property does not have a strong influence on the volume growth
(which can be exponential for instance). In particular, one should
be careful to distinguish it from the well-known {\it doubling
property} stating that there exists a constant $C<\infty$ (not
depending on the radius) such that $V(x,2r)\leq CV(x,r)$ for all
$x\in X$ and $r>0$. Contrary to the locally doubling property, the
doubling property implies polynomial growth, i.e. that there
exists a constant $D<\infty$ such that $V(x,r)\leq r^DV(x,1)$ for
every $x\in X$ and $r\geq 1$.
\end{rem}
For most of the results proved in this paper\footnote{In fact all
the results except the few ones where the infinitesimal structure of
the space is clearly involved. }, we only use the locally doubling
property at scale $r\geq h/2$, if the gradient considered is at
scale $h$. However, to simplify the exposition, we will always
assume that the space is locally doubling.

\subsection*{Examples of locally doubling spaces}

Clearly, the locally doubling property is a very weak property of
controlled geometry:
\begin{itemize}
\item Let $X$ be a connected graph with degree bounded by $d$,
equipped with the counting measure. The volume of balls of radius
$r$ satisfies
$$\forall x\in X,\quad 1\leq V(x,r)\leq d^r.$$ In particular, $X$
is locally doubling. \item Other examples are Riemannian manifolds
with Ricci curvature bounded from below. Assume that the volume of
balls of fixed radius is bounded from above and from below by
constants depending on $r$. Then one can check easily that $X$ is
locally doubling. It is important to note that the locally
doubling property is strictly weaker than this property. One can
easily construct weighted graphs or Riemannian manifolds which are
locally doubling but with unbounded volume for balls of radius
$1$.

\item Let $(X,d,\mu)$ be a metric measure space and let $G$ be a
locally compact group acting by isometries that almost-preserve
the measure, i.e. $$\sup_{g\in G}\sup_{x\in
X}\frac{d(g\cdot\mu)}{d\mu}(x)<\infty.$$ If $G$ acts co-compactly,
then it is easy to check that there exists $C<\infty$ such that
for all $x,y\in X$ and $r>0$,
$$V(x,r)\leq C V(y,r).$$
In particular, $X$ is locally doubling. This obviously applies to
the group itself, equipped with a Haar measure and any metric
which is left-invariant, proper and finite on compact subsets.
\end{itemize}

\begin{defn}\label{largescaledef}
Let $(X,d,\mu)$ and $(X',d',\mu)$ two spaces satisfying the locally
doubling property. Let us say that $X$ and $X'$ are large-scale
equivalent if there is a function $F$ from $X$ to $X'$ with the
following properties
\begin{itemize}

\item[(a)] for every sequence of pairs $(x_n,y_n)\in
(X^2)^{\mathbb{N}}$
$$\left(d(F(x_n),F(y_n))\rightarrow\infty\right)\Leftrightarrow \left(d(x_n,y_n)\rightarrow\infty\right).$$

\item[(b)] $F$ is almost onto, i.e. there exists a constant $C$
such that $[F(X)]_C=X'$.

\item[(c)] For $r>0$ large enough, there is a constant $C_r>0$
such that for all $x\in X$
$$C_r^{-1}V(x,r)\leq V(F(x),r)\leq C_rV(x,r).$$
\end{itemize}
\end{defn}

\begin{crem}
Note that being large-scale equivalent is an equivalence relation
between metric measure spaces with locally doubling property.
\end{crem}

\begin{rem}
If $X$ and $X'$ are quasi-geodesic, then (a) and (b) imply that $F$
is roughly bi-Lipschitz: there exists $C\geq 1$ such that
$$C^{-1}d(x,y)-C\leq d(F(x),F(y))\leq Cd(x,y)+C.$$ This is very
easy and left to the reader. In this case, (a) and (b) correspond to
the classical definition of a {\it quasi-isometry}.
\end{rem}

\begin{ex}
Consider the subclass of metric measure spaces including graphs with
bounded degree, equipped with the countable measure; Riemannian
manifolds with Ricci curvature bounded from below and sectional
curvature bounded from above, equipped with the Riemannian measure.
In this class, quasi-isometries are always large-scale equivalences.
\end{ex}

\section{Examples}\label{Examples}

\subsection{Discretization}\label{ExDiscretization}

Recall that a weighted graph is a connected graph $X$ equipped with
a structure of metric measure space on the set of its vertices, the
distance being the usual geodesic one. Similarly, a weighted
manifold is a Riemannian manifold equipped with a measure $d\mu$
absolutely continuous with respect to the Riemannian measure. A
discretization \cite{Gro',Kan} of a weighted Riemannian manifold $X$
can be defined as a weighted graph large-scale equivalent to $X$.
More generally, a discretization of a metric measure space is a
weighted graph large-scale equivalent to $X$.

Consider some $b>0$. We call a $b$-chain between two points $x,y\in
X$ a chain $x=x_1\ldots x_m=y$ such that for every $1\leq i<m$,
$d(x_i,x_{i+1})\leq b.$ Define another distance on $X$ by setting
$$d_b(x,y)=\inf_{\gamma}l(\gamma)$$
where $\gamma$ runs over every $b$-chains $x=x_0\ldots x_m=y$ and
where $l(\gamma)=\sum_{i=1}^md(x_i,x_{i-1})$ is the length of
$\gamma$.

Let us introduce various natural notions of geodesicity.

\begin{defn}\label{properdef}
We say that a metric space $(X,d)$ is
\begin{itemize}
\item $b$-geodesic if $d(x,y)$ equals the
minimal length of a $b$-chain between $x$ and $y$, or equivalently
if $d=d_b$.

\item quasi-geodesic if there exists $b>0$ such that the identity map $(X,d)\to
(X,d_b)$ is a quasi-isometry;

\item coarsely geodesic if there exists $b>0$ such that
$(X,d_b)\to (X,d)$ is a uniform embedding.
\end{itemize}
Being coarsely geodesic is actually equivalent to being large-scale
uniformly connected (see \cite{tess'}): a space $X$ is large-scale
uniformly connected if there exists $b>0$ such that every $x,y\in X$
can be connected by a $b$-chain whose length only depends on
$d(x,y)$.
\end{defn}

Clearly, being coarsely geodesic is preserved by large-scale
equivalence.

\begin{prop}\label{discretizationprop}
A metric measure space with locally doubling Property admits a
discretization if and only if it is coarsely geodesic. Moreover $X$
is quasi-isometric to a graph if and only if it is quasi-geodesic.
\end{prop}

\noindent{\bf Proof.} Assume that $X=(X,d,\mu)$ is metrically
proper, coarsely geodesic and locally doubling. Consider a minimal
covering of $X$ with balls of radius $h$. We construct a weighted
graph $G(X)$ as follows; the vertices of $G(X)$ are the centers of
the balls; we put an edge between two vertices if the balls
intersect. Since $X$ is coarsely geodesic, $G(X)$ is connected as
soon as $h$ is large enough. Moreover, the coarse geodesicity and
the locally doubling property imply that the injection map
$G(X)\hookrightarrow X$ is a large-scale equivalence. The converse
is obvious. \epr



\subsection{Locally compact groups}\label{ExGroups}

Let $G$ be a group. Recall that a length function on $G$ is function
$L:G\to \R_+$ such that $L(1)=0$ and
$$\forall g,h\in G,\quad L(gh)\leq L(g)+L(h).$$
If $L$ is a length function, then $d(g,h)=L(g^{-1}h)$ defines a
left-invariant pseudo-metric on $G$. Conversely, if $d$ is a
left-invariant pseudo-metric on $G$, then $L(g)=d(1,g)$ defines a
length function on $G$.

\begin{defn} Let $G$ be a locally compact group. A metric
$d$ on $G$ is called uniform if for any of sequence $(g_n,h_n)\in
(G\times G)^{\N}$, $d(g_n,h_n)\to \infty$ if and only if
$g_n^{-1}h_n$ leaves every compact eventually.
\end{defn}
By the Birkhoff-Kakutani metrization theorem \cite[Theorem~7.2]{Hj},
$G$ admits uniform left-invariant metrics if and only if $G$ is
$\sigma$-compact. The following proposition is straightforward and
left to the reader.

\begin{prop} Let $(X,d,\mu)$ and $(Y,d,\mu)$ be metric measure
spaces and let $G$ be a locally compact group acting properly and
co-compactly by isometries that almost-preserve the measure on
both $X$ and $Y$. Then $X$ and $Y$ are locally doubling, and $X$
and $Y$ are large scale equivalent.\;$\blacksquare$
\end{prop}

\begin{cor}
Let $d$ and $d'$ be two uniform metrics on $G$. The spaces $(G,d)$
and $(G,d')$ are doubling at any (large enough) scale and the
identity map $(G,d)\to (G,d')$ is a large scale
equivalence.\;$\blacksquare$
\end{cor}

\begin{defn} Let $G$ be a $\sigma$-compact locally compact group.
The asymptotic class of a metric $d$ is the set of metrics $d'$ on
$G$ such that the identity map $(G,d)\to (G,d')$ is a
quasi-isometry.
\end{defn}

\begin{prop}\label{connectivitegoupprop}
Let $G$ be a locally compact group. The following statements are
equivalent.
\begin{itemize}
\item[(i)] $G$ admits a uniform, coarsely geodesic
metric;

\item[(ii)] $G$ admits a uniform, quasi-geodesic metric;

\item[(iii)] $G$ admits a left-invariant, proper, quasi-geodesic metric;

\item[(iv)] $G$ admits a left-invariant proper metric, quasi-isometric to a graph with bounded degree;

\item[(v)] $G$ is compactly generated.
\end{itemize}
\end{prop}
\bpr Clearly, $(iii)\Rightarrow(ii)\Rightarrow(i)$ are obvious,
$(iii)\Leftrightarrow(iv)$ results from
Proposition~\ref{discretizationprop}. Let us prove that
$(v)\Rightarrow (iv)$. Assume that $G$ is compactly generated and
let $S$ be a compact symmetric subset $S$. One can equip $G$ with a
uniform quasi-geodesic length function setting
$$\forall g\in G,\quad |g|_S=\inf\{n\in \N, g\in S^n\}.$$
Now, let us prove that $(i)\Rightarrow(v)$. Suppose that $G$ has a
uniform, coarsely geodesic metric $d$ with constant $C$. Since $d$
is uniform, there exists $R<\infty$ such that for all $g\in G$, the
closed ball $B(g,C)$ is compact and contained in $g\cdot B(1,R).$

We claim that $G$ is generated by $B(1,R)$. Fix $g\in G$. Indeed,
let $g_1=1, \ldots, g_n=g$ be a chain such that $d(g_i,g_{i+1})\leq
C$ for every $1~\leq i\leq n-1$. We have $g_{i+1}\in B(g_i,C)\subset
g_i\cdot B(1,R)$. Hence, an immediate induction shows that $g\in
B(1,R)^n$ and we are done.\epr

\section{Equivalence of Sobolev inequalities with respect to different
gradients}\label{equivalentgradientsection}

Here, we show that {\it large-scale} Sobolev inequalities do not
really depend on the kind of gradient that we use to write them. In
spite of its easy and short proof, this result is crucial for our
purpose since it shows that our definitions are natural.

The following proposition results immediately from the definitions.
\begin{prop}\label{comparaisonfacile}
If $h'\geq h>0$, then $$\||\nabla f|_{h'} \|_p\geq \||\nabla f|_h
\|_p .$$ Moreover, if $P$ is a viewpoint at scale $h$ with constants
$c$ and $A$ (see the definition below) and if $q\leq q'\leq\infty$,
then
$$c|\nabla f|_{h,q}\leq |\nabla f|_{P,q}\leq |\nabla f|_{P,q'}\leq |\nabla f|_{Ah}\quad \forall f\in L^{\infty}(X).~\blacksquare$$
\end{prop}

The non-trivial comparisons between different gradient are
summarized in the following proposition.

\begin{prop}\label{comparaison}
Let $X$ be some metric measure space satisfying a Sobolev inequality
$(S_{\varphi}^p)$ at scale $h$. Then, for any viewpoint
$P=(P_x)_{x\in X}$ at scale $2h$, $X$ satisfies $(S_{\varphi}^p)$
w.r.t. $|\nabla|_{P,q}$ for any $q\geq 1.$
\end{prop}

\bpr By Proposition \ref{comparaisonfacile}, it suffices to prove
that $X$ satisfies $(S_{\varphi}^p)$ w.r.t. $|\nabla|_{2h,1}$. Write
$$P_x=\frac{1}{V(x,h)}1_{B(x,h)}\quad \forall x\in X.$$ For every
$f\in  L^{\infty}(X)$ we write
$$Pf(x)=\int fdP_x, \quad \forall x\in X.$$

\begin{lem}\label{compargradfgradPf}
There exists $C<\infty$ such that
$$|\nabla Pf|_h(x)\leq C|\nabla f|_{h,1}(x)\quad \forall f\in L^{\infty}(X), \forall x\in X.$$
\end{lem}
\bpr Consider some $y\in B(x,h).$
$$|Pf(x)-Pf(y)|\leq |Pf(x)-f(x)|+|Pf(y)-f(x)|\leq C|\nabla|_{2h,1}f(x).$$
with $C<\infty$ depending only on the doubling constant at scale
$h$.\epr

\

\noindent Now apply the Sobolev inequality $(S_{\varphi}^p)$ at
scale $h$ to $Pf$,
$$\||\nabla Pf|_h\|_p\geq \varphi^{-1}(\Omega)\|Pf\|_p\geq  \varphi^{-1}(\Omega)\|f\|_p- \varphi^{-1}(\Omega)|\|f\|_p-\|Pf\|_p|.$$
Now, if $\||\nabla f|_{h,1}\|_p\geq \|f\|_p/2$, there is nothing to
prove. Hence, assuming the contrary, and since
$|\|f\|_p-\|Pf\|_p|\leq \||\nabla f|_{h,1}\|_p$, we obtain
$$\||\nabla Pf|_h\|_p\geq \varphi^{-1}(\Omega)\|f\|_p/2,$$
which yields
$$\||\nabla f|_{h,1}\|_p\geq C^{-1}\varphi^{-1}(\Omega)\|f\|_p/2$$
thanks to the lemma.\epr

\section{Invariance of Sobolev inequalities under large-scale equivalence}

The aim of this section is to prove the following theorem.

\begin{thm}\label{largescalethm}
Let $F:X\rightarrow X'$ be a large-scale equivalence between two
spaces $X$ and $X'$ satisfying the locally doubling property. Assume
that for $h>0$ fixed, the space $X$ satisfies a Sobolev inequality
$(S^p_{\varphi})$ at scale $h$, then there exists $h'$, only
depending on $h$ and on the constants of $F$ such that $X'$
satisfies $(S_{\varphi}^p)$ at scale $h'$. In particular,
large-scale Sobolev inequalities are invariant under large scale
equivalence.
\end{thm}

To prove Theorem~\ref{largescalethm}, we will first prove some
preliminary results.
\subsection{Thick subsets}

\begin{defn}
A subset $A$ of a metric space is called $h$-thick if it is a
reunion of closed balls of radius $h$.
\end{defn}

Roughly speaking, the following proposition says that Large-scale
Sovolev inequalities hold if and only if they hold for functions
with thick support.

\begin{prop}\label{thickprop}
Let $X=(X,d,\mu)$ be a metric measure space. Fix some $h>0$ and some
$p\in [1,\infty]$. There exists a constant $C>0$ such that for any
$f\in   L^{\infty}(X)$, there is a function $\tilde{f}\in
L^{\infty}(X)$ whose support is included in a $h/2$-thick subset
$\Omega$ such that
$$\mu(\Omega)\leq \mu(Supp(f))+C$$
and for every $p\in [1,\infty]$,
$$\frac{\| |\nabla \tilde{f}|_{h/2} \|_p}{\|\tilde{f}\|_p}\leq C\frac{\| |\nabla f|_{h}\|_p}{\| f\|_p}.$$
\end{prop}

\bpr Let us prove the proposition for $p<\infty$. Let $f\in
  L^{\infty}(X)$ be such that $\| f\|_p=1$. Assume that $f$ satisfies
$$\| |\nabla|_h f\|_p\geq \frac{1}{2}.$$
Then, for $\tilde{f}$, consider for instance the indicator function
of a ball $B(x,a)$ of volume $1$ (so that $\| \tilde f\|_p=1)$. We
have
$$\| |\nabla\tilde{f}|_{h/2} \|_p^p\leq \mu(B(1+h/2))\leq C\mu(B(x,a))=C.$$
Thus, let us assume that
$$\| |\nabla f|_h\|_p\leq \frac{1}{2}.$$
Let $\Omega$ be the subset of $Supp(f)$ defined by
$$\Omega=\{x\in X, d(x,Supp(f)^c)\geq h/2\}$$
and set
$$\tilde{f}=f\cdot1_{\Omega}.$$
Note that for every $x\in Supp(f)\smallsetminus \Omega$, there
exists some $y\in B(x,h)$ such that $f(y)=0$. Therefore, we have
$|f(x)|\leq |\nabla f|_h(x).$ Hence,
$$\int_X|\tilde{f}|^pd\mu\geq \int_X|f|^pd\mu -\int_{X} (|\nabla f|_h)^pd\mu\geq \frac{1}{2}.$$
On the other hand, let $x\in \Omega$. If $d(x,Supp(f))\geq h$, then
$$|\nabla\tilde{f}|_{h/2}= |\nabla|_{h/2}f\leq |\nabla|_h f.$$
Otherwise,
$$|\nabla\tilde{f}|_{h/2}\leq \max\left\{|f(x)|,\sup_{y\in B(x,h/2)}|f(x)-f(y)|\right\}$$
and
$$|\nabla f|_{h}= \sup_{y\in B(x,h)}|f(x)-f(y)|=\max\left\{|f(x)|,\sup_{y\in B(x,h)}|f(x)-f(y)|\right\}.$$
Thus
$$|\nabla\tilde{f}|_{h/2}\leq |\nabla f|_{h};$$
so we are done.\epr


\

On the other hand, the locally doubling property ``extends" to thick
subsets in the following sense.

\begin{prop}\label{DVthickprop}
Let $X$ be a metric measure space satisfying the locally doubling
property. Fix two positive numbers $u$ and $v$. There exists a
constant $C=C(u,v)<\infty$ such that for any $u$-thick subset
$A\subset X$, we have
$$\mu([A]_v)\leq C\mu(A).$$
\end{prop}
\bpr The proof follows from standard covering arguments.\epr

\subsection{Rough volume-preserving property}

Let us prove a useful {\it rough} volume preserving property of
large scale equivalences.

\begin{prop}\label{volpres}
Let $X=(X,d,\mu)$ and $X'=(X',d',\mu')$ be two spaces satisfying the
locally doubling property and let $F:X\rightarrow X'$ be a
large-scale equivalence. Let $u>0$, then there exists a constant
$C=C(u,F)$ such that

\noindent (1) If $A\subset X$ and $A'\subset X'$ are such that
$[F^{-1}(A')]_u\subset A$, then $\mu'(A')\leq C\mu(A).$

\noindent (2) If $A\subset X$ and $A'\subset X'$ are such that
$[F(A)]_u\subset A'$, then $\mu(A)\leq C\mu'(A').$

\end{prop}
\bpr Let us prove (1). Let $Z$ be a maximal set of $2u$-separated
points of $F^{-1}(A')$. Clearly, the balls $(B(z,u))_{z\in Z}$ are
disjoint and included in $A$. On the other hand, maximality of $Z$
implies that the family $(B(z,2u))_{z\in Z}$ forms a covering of
$A$. So we have
\begin{equation}\label{eq1}
\sum_{z\in Z}\mu(B(z,u))\leq \mu(A)\leq \sum_{z\in Z}\mu(B(z,2u))
\end{equation}
By property $(a)$ of a large-scale equivalence, there exists $v$
such that for every $x\in X$, $F(B(x,2u))\subset B(F(x),v).$ In
particular, the family $((B(F(z),v))_{z\in Z}$ forms a covering of
$F(A)$. Using Property (c) of a large-scale equivalence and Doubling
Property at any scale of $X$ together with (\ref{eq1}), we get
\begin{eqnarray*}
\mu(A')\leq \mu'(F(A))\leq \sum_{z\in Z}\mu'(B(F(z),v))&\leq &
C'\sum_{z\in Z}\mu(B(z,v))\\ & \leq & C\sum_{z\in Z}\mu(B(z,u))\leq
C\mu(A)
\end{eqnarray*}
which proves the proposition. \epr

\subsection{Proof of the invariance under large-scale
equivalence}\label{largescaleinvariancesection}

Let $F:X\rightarrow X'$ be a large-scale equivalence between two
spaces $X$ and $X'$ satisfying the locally doubling property. Assume
that $f\in  L^{\infty}(X')$. For every $h>0$, define a function on
$X$
$$\forall x\in X, \quad \psi_h(x)=\sup_{y\in B(x,h)}|f\circ F(y)|.$$

\begin{lem}\label{lem1}
For $h$ large enough, there exists a constant $c=c(h,f)>0$ such that
$$\mu(\{\psi_h^p\geq t\})\geq c\mu'(\{|f|^p\geq t\}).$$ In
particular, for every $p>0$,
$$ \| \psi_h\|_p\geq c\| f\|_p.$$
\end{lem}
\bpr We can obviously assume that $p=1$ and that $f\geq 0$. Thanks
to Proposition \ref{volpres}, we only have to check that
$$[F^{-1}(\{f\geq t\})]_{h}\subset\{\psi_h\geq t\}.$$
Indeed, let $x\in F^{-1}(\{f\geq t\})$. Then $f\circ F(x)\geq t$. So
for all $y\in B(x,h)$, we have $\psi_h(y)\geq t$.\epr

\begin{lem}\label{lem2}
For $h'$ large enough, there exists a constant $C<\infty$ such that
$$\mu(\{(|\nabla\psi_h|_h)^q> t\})\leq C\mu'(\{|(\nabla f|_{h'})^q>
t/2\}).$$ In particular, for every $q>0$,
$$\||\nabla\psi_h|_h\|_q\leq C\|
|\nabla f|_{h'}\|_q.$$

\end{lem}
\bpr We can of course assume that $q=1$. Thanks to Proposition
\ref{volpres}, it suffices to prove that for $h'$ large enough,
$$[F(\{|\nabla\psi_h|_h> t\})]_{h'/2}\subset\{|\nabla f|_{h'}>t/2\}.$$
Indeed, let $x\in X$ be such that $|\nabla\psi_h|_h(x) > t$. This
means that there exists $y\in B(x,h)$ such that $|f\circ F(x)-f\circ
F(y)|>t$. On the other hand, by property (a) of a large-scale
equivalence, one can choose $h'$ such that $d(F(x),F(y))\leq h'/2$.
Hence,
$$\forall z\in B(F(x),h'/2),\quad |\nabla f|_{h'}(z)\geq
\max\{|f(x)-z|,|f(y)-z|\}\geq t/2.$$ So $z\in \{|\nabla
f|_{h'}>t/2\}$. \epr
\begin{lem}\label{lem3}
For $u$ large enough, there exists a constant $C<\infty$ such that
$$\mu'\left(Supp(\psi_h)\right)\leq C\mu\left([Supp(f)]_u\right).$$
\end{lem}
\bpr This follows trivially from Proposition \ref{volpres}.

\

\noindent{\bf Proof of Theorem~\ref{largescalethm}} Let $\Omega$ be
a compact subset of $X'$ of measure $m$. We want to prove that every
$f\in L^{\infty}(\Omega)$ satisfies
$$\| f\|_p\leq C\varphi(Cm)\| |\nabla f|_{h}\|_p$$ with $h'$ and $C$ depending only on $F$, $h$ and $X$.
Thanks to Proposition \ref{thickprop} and up to choose a larger
$h'$, we can assume that $\Omega$ is $v$-thick for any $v>0$. Then,
thanks to Lemma \ref{lem3} and to Proposition \ref{DVthickprop}, we
have
$$Supp(\psi_h)\leq C'm$$
for some constant $C'$. So apply $(S_{\varphi}^p)$ to $\psi_h$ and
then conclude thanks to Lemmas \ref{lem1} and \ref{lem2}.\epr

\section{Sobolev inequality
$(S_{\varphi}^2)$  and on-diagonal upper bounds for random
walks}\label{random}

In this section, we revisit the relations (see \cite{C} for a
survey) between Sobolev inequalities $(S_{\varphi}^2)$ and
on-diagonal upper bounds for random walks in our general context.
The main purpose is to prove a version of
\cite[Theorem~7.2]{C}(see also \cite[Theorem~7.2]{Coulhon}) to our
more general context.

\begin{thm}\label{random/sobolevThm}
Let $X=(X,d,\mu)$ be a metric measure space and let $P=(P_x)_{x\in
X}$ be a symmetric view-point at scale $h$ on $X$. Let $\varphi$ be
some increasing positive function. Define $\gamma$ by
$$t=\int_{0}^{1/\gamma(t)}\left(\varphi(v)\right)^2\frac{dv}{v}.$$
\begin{itemize}
\item[(i)] Assume that $X$ satisfies a Sobolev inequality
$(S_{\varphi}^2)$ w.r.t. $|\nabla f|_{P^2,2}$. Then
$$p_x^{2n}(x)\leq \gamma(cn)\quad  \forall n\in \N,$$
for some constant $c>0$.

\item[(ii)] If the logarithmic derivative of $\gamma$ has at most polynomial growth and if
$$p_x^{2n}(x)\leq \gamma(n)\quad \forall n\in \N,$$ then $X$ satisfies
$(S_{\varphi}^2)$ w.r.t. $|\nabla|_{P,2}.$
\end{itemize}
\end{thm}

Similarly we have the following version\footnote{The proofs are
straightforward adaptations of their versions for graphs in \cite{C}
so we will only prove Theorem~\ref{random/sobolevThm}.} of
\cite[Theorem~7.1]{C}

\begin{thm}\label{{random/sobolevThm}}
Let $X=(X,d,\mu)$ be a metric measure space and let $P=(P_x)_{x\in
X}$ be a symmetric view-point at scale $h$ on $X$. Define $\gamma$
by
$$t=\int_{0}^{1/\gamma(t)}\left(j_{X,2}(v)\right)^2\frac{dv}{v}.$$
where $j_{X,2}$ is the isoperimetric profile of $X$ defined with the
gradient $|\nabla f|_{P^2,2}$. If the logarithmic derivative of
$\gamma$ has at most polynomial growth, then there exists a constant
$C>0$ such that
$$\gamma(Cn)\leq \sup_{x\in X}p_x^{2n}(x)\leq \gamma(n)\quad  \forall n\in \N.$$
\end{thm}

\noindent{\bf Proof of Theorem~\ref{random/sobolevThm}.} In
\cite[Theorem~7.2]{C}, the same result is proved for a weighted
graph $(X,\mu)$ using the usual notion of gradient on graphs (see
Example~\ref{graphEx}) and where $P$ is the standard random walk on
$(X,\mu)$. Their proof only relies on the following formal link
between $P$ and the gradient.
$$ c(\|f\|_2^2-\|Pf\|_2^2)\leq \|\nabla f\|_2^2\leq C(\|f\|_2^2-\|Pf\|_2^2).$$
Here, this relation is satisfied when considering the gradient
$|\nabla|_{P^2,2}$ and we even have the equality
\begin{lem}
For every $f\in L^2(X)$, we have $$\||\nabla
f|_{P^2,2}\|_2^2=\|f\|_2^2-\|Pf\|_2^2.$$
\end{lem}
\bpr We have (see section~\ref{laplacianSection})
\begin{eqnarray*}
\||\nabla f|_{P^2,2}\|_2^2& = & \langle \Delta_{P^2}f,f\rangle \\
                          & = &\langle (id-P^2)f,f\rangle\\
                          & = & \|f\|_2^2-\langle P^2f,f\rangle\\
                          & = & \|f\|_2^2-\langle Pf,Pf\rangle \; \blacksquare
\end{eqnarray*}

\

So the proof of \cite[Theorem~7.2]{C} can be used formally in our
context. However, for the sake of completeness, we give a sketch of
this proof. First, using that $P^{n}$ is symmetric, one checks
easily that
$$\sup_{x\in X}p_x^{2n}(x)=\|P^{2n}\|_{1\to \infty}$$
where $\|\cdot\|_{p\to q}$ denotes the operator norm form
$L^p(X,\mu)$ to $L^q(X,\mu)$.

\noindent{\bf Proof of (i).} Assume that $(S_{\varphi}^2)$ holds.
Let us start with an important lemma.

\begin{lem}\label{Nashlem}
The Sobolev inequality $(S_{\varphi}^2)$ for the $L^2$-gradient
w.r.t. the viewpoint $P$ is equivalent to the so-called Nash
inequality
$$
\|f\|_2^2\leq
C\varphi^2\left(C\frac{\|f\|_1^2}{\|f\|_2^2}\right)\||\nabla|_{P,2}\|_2^2.$$
\end{lem}
\bpr Assume that a function $f$ satisfies Nash inequality.  Using
Schwarz inequality and the fact that $\varphi$ is nondecreasing, we
obtain
$$\|f\|_2^2\leq \varphi^2\left(\frac{\|f\|_2^2}{\|f\|_1^2}\right)\||\nabla f|_{P^2,2}\|_2^2\leq \varphi^2(|\Omega|)\||\nabla f|_{P^2,2}\|_2^2.$$
The proof of the other implication relies on an argument of
Grigor'yan in \cite{Gri}. Assume that $(S_{\varphi}^2)$ holds. Let
$f\in   L^{\infty}(X)$. For every $\lambda>0,$ since
$f<2(f-\lambda)$ on $\{f>2\lambda\}$, we may write
\begin{eqnarray*}
\int f^2 & = & \int_{f>2\lambda} f^2+\int_{f\leq 2\lambda} f^2\\
         & \leq & 4\int_{f>2\lambda} (f-\lambda)^2+2\lambda\int_{f\leq 2\lambda}
         f\\
         & \leq & 4\int_{f>2\lambda} (f-\lambda)^2+2\lambda\|f\|_1
\end{eqnarray*}
Now applying $(S_{\varphi}^2)$ to $(f-\lambda)_+$ gives
$$\int(f-\lambda)_+^2 \leq \varphi^2(\mu(\{f>2\lambda\}))\||\nabla
f|_{P^2,2}\|_2^2, $$ that is, since
$$\mu(\{f>\lambda\})\leq \frac{\|f\|_1}{\lambda}$$
and $\varphi$ is non-decreasing,
$$\int (f-\lambda)_+^2\leq \varphi^2\left(\frac{\|f\|_1}{\lambda}\right)\||\nabla f|_{P^2,2}\|_2^2.$$
Therefore
$$\int f^2\leq 4\varphi^2\left(\frac{\|f\|_1}{\lambda}\right)\||\nabla f|_{P^2,2}\|_2^2+2\lambda\|f\|_1.$$
Letting $\eps>0$ and taking $\lambda=\eps\|f\|_2^2/\|f\|_1$ in this
equation yields
$$\|f\|_2^2\leq 4\varphi^2\left(\frac{\|f\|_2^2}{\eps\|f\|_1^2}\right)\||\nabla f|_{P^2,2}\|_2^2+2\eps
\|f\|_2^2$$ or equivalently,
$$\|f\|_2^2\leq \frac{4}{1-2\eps}\varphi^2\left(\frac{\|f\|_2^2}{\eps\|f\|_1^2}\right)\||\nabla f|_{P^2,2}\|_2^2$$
Taking $\eps=1/4$, for example yields
$$\|f\|_2^2\leq 8\varphi^2\left(4\frac{\|f\|_2^2}{\|f\|_1^2}\right)\||\nabla f|_{P^2,2}\|_2^2$$
which is the expected Nash inequality. \epr

\

Now, consider $f\in L^1(X,\mu)$, non-negative, with $\|f\|=1$ and
define a sequence $u_n=\|P^{n}f\|_2^2$. The above inequality applied
to the function $P^{n}f$ thus reads as
$$u_n\leq \varphi^2(1/u_n)(u_{n}-u_{n+1})$$
since $\|P^{n}f\|_1=\|f\|_1=1$ by Markov property of $P$. Let $t\to
u_t$ be the increasing, piecewise linear function extending $u_n$ on
$\R_+$. If we put $v_t=1/u_t$, then the above inequality becomes
$$dt\leq \varphi^2(v_t)\frac{dv_t}{v_t},$$ hence, by integrating
between $0$ and $t$, we obtain
$$t\leq \int_{v_0}^{1/v_t}\varphi^2(s)\frac{ds}{s};$$
and since by definition
$$t=\int_{0}^{1/\gamma(t)}\left(\varphi(v)\right)^2\frac{dv}{v},$$
this means that $\gamma(t)\leq v_t$, i.e.
$$\|P^{n}f\|_2^2\leq \gamma(n)$$
from which we deduce
$$\|P^n\|_{1\to 2}\leq \sqrt{\gamma(n)}.$$
Now, using the fact that $P^n$ is symmetric,
$$\|P^{n}\|_{2\to\infty}=\|P^n\|_{1\to 2}\leq \sqrt{\gamma(n)}.$$
Hence
$$\|P^{2n}\|_{1\to \infty}\leq \|P^n\|_{2\to\infty}\|P^n\|_{1\to
2}\leq \gamma(n).$$ So (i) follows.

\noindent{\bf Proof of (ii).} Assume that the decay
$\|P^{2n}\|_{1\to \infty}\leq \gamma(n)$ holds. Observe that
$\|P^{2n}\|_{1\to \infty}=\|P^{n}\|_{1\to 2}$, then take $f$ with
$\|f\|_1=1$ and define as above $u_n=\|P^nf\|_2^2.$ Since $P$ is
self-adjoint,
$$\|P^nf\|_2^2=\langle P^nf,P^nf\rangle=\langle P^{n-1}f,P^{n+1}f\rangle\leq \|P^{n-1}f\|_2\|P^{n+1}f\|_2.$$
In other words, $u_n^2\leq u_{n-1}u_{n+1}$ and $u_{n+1}/u_n$ is
nondecreasing in $n$. It follows that
$$\left(\frac{u_1}{u_0}\right)^n\leq \frac{u_1}{u_0}\frac{u_2}{u_1}\ldots \frac{u_n}{u_{n-1}}=\frac{u_n}{u_0}.$$
Now, since by assumption $u_n\leq \gamma(n)$,
$$\log\frac{\|f\|_2^2}{\gamma(n)}\leq \log\frac{u_0}{u_n}\leq n\log\frac{u_0}{u_n}\leq n\left(\frac{u_0}{u_1}-1\right),$$
hence
$$\|Pf\|_2^2\leq \left(\frac{n}{\log\frac{\|f\|_2^2}{\gamma(n)}}\right)(\|f\|_2^2-\|Pf\|_2^2), \; \forall n\in \N.$$
Finally, for all $f$ such that $\|f\|_1=1$,
$$\|f\|_2^2\leq \left(\frac{n}{\log\frac{\|f\|_2^2}{\gamma(n)}}+1\right)(\|f\|_2^2-\|Pf\|_2^2), \; \forall n\in \N.$$
An optimization\footnote{This is where condition $(\delta)$ is
needed.} in $n$ yields the Nash inequality that is equivalent to
$(S^{2}_{\varphi}$) by Lemma~\ref{Nashlem}. \epr

\section{Controlling the scale of Sobolev Inequalities}\label{controlscalesection}

\subsection{Going down the scale}

In this section, we address the following question. Let
$X=(X,d,\mu)$ be a metric measure space $X$ satisfying a Sobolev
inequality at scale $h$; we know that it automatically satisfies the
same Sobolev inequalities at any larger scale; but under what
assumptions does it satisfy this inequality at some smaller scale
$h'$? This can be compared to a similar discussion in \cite{tess'}
where we considered the isoperimetric properties of a metric measure
space\footnote{This is a particular case of the present discussion
corresponding to $p=1$.}.

For example, consider $X=\mathbb{Z}^d$ ($d\geq 2$) equipped with the
distance $d(x,y)=\sum_{i=1}^d |y_i-x_i|$ and with the countable
measure. It is well known that $X$ satisfies a Sobolev inequality
$S(d/(d-1),1)$ at any scale $\geq 1$. But no Sobolev inequality is
available at a scale $s<1$ since for every $f\in L^{\infty}(X)$,
$|\nabla f|_{s}=0$. Clearly, the problem comes from the lack of
connectivity at scale $<1$.

The following proposition shows that Property\footnote{see
definition~\ref{properdef}.} of coarse $b$-geodesicity (also called
uniform $b$-connectedness) together with Property of locally
doubling are sufficient to control the minimal scale at which
Sobolev inequalities may be valid.

\begin{prop}\label{scale}
Assume that $X$ is a  coarse $b$-geodesic space satisfying the
locally doubling property $r\geq b$. Then $X$ satisfies a
large-scale Sobolev inequality if and only if it satisfies the
same Sobolev inequality at scale $2b$ (but with different
constants). In other words, the asymptotic behavior of the
isoperimetric profile $j_{X,p}$ does not depend on the scale,
provided it is larger than $2b$.
\end{prop}

\bpr Let $f\in L^{\infty}(X)$. Let us prove that for all $h\geq
2b$, there is a constant $C=C(h)<\infty$ such that for every $t>0$
\begin{equation}\label{eq2}
\mu(\{|\nabla f|_h> t\})\leq C\mu(\{|\nabla f|_{2b}>t/C\}).
\end{equation}
Consider a point $x\in \{|\nabla f|_{h}> t\}$: there is $y\in
B(x,h)$ such that $|f(x)-\varphi(y)| >t$. Now, let $x=x_1\ldots
x_m=y$ be a $b$-connecting chain between $x$ and $y$ (with $m$ only
depending on $h$). Clearly, there exists $1\leq i<m$ such that
$|\varphi(x_i)-\varphi(x_{i+1})|> t/m$. So in particular, for all
$z\in B(x_i,b)$, $|\nabla f|_{2b}(z)>t/(2m)$. Let $Z$ be a maximal
$2E$-separated subset of $\{|\nabla f|_{h}> t\}$. The balls
$(B(z,2E))_{z\in Z}$ form a covering of $\{|\nabla f|_{h}> t\}$. On
the other hand, by the previous discussion, in each ball $B(z,E)$,
one can find a ball $B(x_z,b)$ included in $\{|\nabla f|_{2b}>
t/(2m)\}$. Since the balls $(B(x_z,b))_{z\in Z}$ are disjoint,
(\ref{eq2}) follows from locally doubling property $r\geq b$. \epr

\medskip

As an interesting corollary of Proposition~\ref{scale} and
Proposition~\ref{connectivitegoupprop}, we obtain that if $h$ is
large enough, a Sobolev inequality is satisfied at scale $h$ on a
locally compact {\it compactly generated} group if and only if it is
satisfied at large scale. It also allows to define an
$L^p$-isoperimetric profile on locally compact compactly generated
groups, whose asymptotic behavior does not depend on the scale,
provided it is large enough. As a corollary of
Theorem~\ref{largescalethm}, we therefore have

\begin{cor}\label{groupisocor1}
Let $H$ and $G$ be quasi-isometric amenable unimodular locally
compact compactly generated group. Then,

\noindent(1) $ j_{H,p}\approx j_{G,p}$;

\noindent(2) $J_{H,p}^b\approx J_{G,p}^b.~\blacksquare$
\end{cor}

\begin{rem}
In particular, for $p=2$, the asymptotic behavior when $r\to
\infty$ of first eigenvalue $\delta_{P}(r)$ of the Laplacian
associated to any viewpoint $P$ at scale $1$, acting on
square-integrable functions supported in balls of radius $r$, does
not depend on $P$. We therefore denote it by $\delta_G$. Part (2)
of Corollary~\ref{groupisocor} for $p=2$ says that the asymptotic
behavior of $\delta_G$ is invariant under quasi-isometry (see
Section~\ref{laplacianSection}).
\end{rem}

\subsection{From finite scale to infinitesimal scale}\label{finite->infinitesimalsection}

\begin{defn}\label{uppergradDef}(see for instance
\cite[Definition~1.18]{Sem}) Let $(X,d)$ be a metric space, and let
$u$ and $g$ be two Borel measurable functions defined on $X$, with
$u$ real-valued and $g$ taking values in $[0,\infty]$. We say that
$g$ is an {\it generalized gradient} of $u$ if
$$|u(\gamma(a))-u(\gamma(b))|\leq \int_a^bg(\gamma(t))dt$$
whenever $a,b\in \R$ and $\gamma: [a,b]\to X$ is $1$-Lipschitz (so
that $d(\gamma(s),\gamma(t))\leq |s-t|$ for all $s,t\in [a,b]$).
\end{defn}

\begin{ex}\label{standarduppergrad}\cite[Lemma~1.20]{Sem} The function $g$ defined by
$$g(x)=\liminf_{r\to 0}r^{-1}\sup_{y\in B(x,r)}|u(y)-u(x)|$$
is a generalized gradient of $u$. Let us call $g$ the standard upper
gradient of $u$ and we denote it by $|\overline{\nabla}u|$.
\end{ex}

The following proposition is obvious by passing to the limit.

\begin{prop}
Fix $p\in[1,\infty]$. Assume that for every $h>0$, $(X,d,\mu)$
satisfies a Sobolev inequality $(S_{\varphi}^p)$ w.r.t. the gradient
$|\nabla|_{h}$. Suppose that the constants appearing in these
inequalities are uniform with respect to $h$, then $X$ satisfies
$(S_{\varphi}^p)$ w.r.t. the standard upper gradient.
\end{prop}

The following fact had already been noticed in the case of a
discretization of a manifold \cite{CS}. Its proof, here, is
straightforward from the definition of $|\nabla|_{P,p}$.

\begin{prop}\label{finite->infinitesimal}
Fix some $h>0$  and $p\in[1,\infty]$. Let $(X,d,\mu)$ be a metric
measure space with doubling property at radius $\geq h$, and let $P$
be a viewpoint at scale $h$ on $X$. Suppose that a function $u\in
L^{\infty}(X)$ satisfies $(S_{\varphi}^p)$ w.r.t. $|\nabla|_{P,p}$.
Let $g$ be an generalized gradient of $u$. We assume that $u$
satisfies the following local Poincaré inequality $(P(1,p))_{loc}$
$$\int_{B(x,h)}|h(y)-h(x)|^pdP_x(y)\leq C\int_{B(x,h')}g^p(y)d\mu(y)$$
for some constants $C,h<\infty$. Then $u$ satisfies
$(S_{\varphi}^p)$ w.r.t. $g$.
\end{prop}

\begin{ex}
Let $M$ be a Riemannian manifold. Then the local norm of its usual
gradient trivially coincides with the standard upper gradient on
$M$. Now, assume that $M$ satisfying a local Poincaré inequality (as
in the Proposition) and let $X$ be a discretization of $M$.
According to Theorem~\ref{largescalethm}, if $X$ satisfies
$(S_{\varphi}^p)$, then $M$ also satisfies $(S_{\varphi}^p)$ w.r.t.
its usual gradient.
\end{ex}

\subsection{From infinitesimal scale to finite scale}\label{infinitesimal->finitesection}

In this last section, we will prove that if a metric measure space
satisfies a Sobolev inequality w.r.t. the standard upper gradient
(see Exemple~\ref{standarduppergrad}), then it satisfies this
Sobolev inequality at any scale.

\begin{thm}\label{infinitesimal->fini}
Fix $p\in [1,\infty].$ Let $(X,d,\mu)$ be a metric measure space
satisfying the locally doubling property. Assume that $(X,d,\mu)$
satisfies a Sobolev inequality $(S_{\varphi}^p)$ w.r.t. the standard
upper gradient $|\overline{\nabla}|$. Then $X$ satisfies
$(S_{\varphi}^p)$ w.r.t. $|\nabla|_{h}$ for every $h>0$.
\end{thm}
\bpr Assume that $X$ satisfies  $(S_{\varphi}^p)$ w.r.t. the
standard upper gradient. Using the same tools as in the proof of
Proposition~\ref{comparaison}, one can see that it suffices to show
that for every $h>0$ and every function $f$, there exists a
viewpoint $P$ at scale $h/2$ such that
\begin{equation}\label{Sobofaible}
\|Pf\|_p\leq C\varphi(\mu(\Omega))\||\nabla Pf|_{h}\|_p
\end{equation}
where $\Omega$ is a measurable subset containing the support of
$f$. According to Proposition~\ref{thickprop}, we can assume that
$Supp(f)$ is thick. Thus, thanks to Proposition~\ref{DVthickprop},
we can replace $\Omega$ by $[\Omega]_{Ah}$ that\footnote{A is the
large constant appearing in the definition of a viewpoint at scale
$h$.} contains $Supp(Pf)$ for any viewpoint $P$ at scale $h/2$
Hence, it suffices to prove that $(S_{\varphi}^p)$ w.r.t.
$|\nabla|_h$ is satisfied for functions of the form $Pf$, with
$f\in L^{\infty}(X)$.

Define a $1$-Lipschitz map $\theta: X\times X\to \R_+$ by
$\theta(x,y)=d(y,B(x,h)^c)$. Write
$$p_x(y)=\frac{\theta(x,y)}{K(x)},$$where
$K(x)=\int_{B(x,h)}\theta(x,z)d\mu(z)$. Since $X$ is locally
doubling, one can easily check that $p_x(y)$ is the density of a
viewpoint $P$ at scale $h$. Moreover, $D^{-1}V(x,h)\leq K(x')\leq
DV(x,h)$ where $D\geq 1$ only depends on the doubling constant at
scale $h$.

Let $x'$ be a point distinct from $x$. We have
\begin{eqnarray*}
Pf(x')-Pf(x) & =
&\int_{X}(p_{x'}(y)-p_x(y))f(y)d\mu(y)\\ & = & \int_{X}(p_{x'}(y)-p_x(y))(f(y)-f(x))d\mu(y)\\
& = & \int_X
\frac{\theta(x',y)K(x)-\theta(x,y)K(x')}{K(x)K(x')}(f(y)-f(x))d\mu(y)\\
& = & \int_X \frac{(\theta(x',y)-\theta(x,y))K(x)-\theta(x,y)\left(K(x')-K(x)\right)}{K(x)K(x')}(f(y)-f(x))d\mu(y)\\
\end{eqnarray*}
Since $X$ is locally doubling, it is not difficult to see that for
$x'$ closed to $x$, $C^{-1}K(x)\leq K(x')\leq CK(x)$ where $C\geq 1$
only depends on the doubling constant at scale $h$. Hence,
\begin{eqnarray*}
|\overline{\nabla}Pf|(x) & \leq &
C\int_X \frac{|\overline{\nabla}_x\theta|(x,y)K(x)+\theta(x,y)|\overline{\nabla}K|(x)}{K(x)^2}|f(y)-f(x)|d\mu(y)\\
\end{eqnarray*}
On the other hand, note that
$$|\overline{\nabla}K|(x)\leq \int_X |\overline{\nabla}_x\theta|(x,z)d\mu(z)\leq V(x,h).$$
Up to change the constant $C$, we conclude that
\begin{eqnarray*}
|\overline{\nabla}Pf|(x) & \leq & C
\frac{1}{V(x,h)}\int_{B(x,h)}|f(y)-f(x)|d\mu(y)\\
                         & \leq C|\nabla f|_{h}(x).
\end{eqnarray*}

Now, to conclude, it remains to apply $(S_{\varphi}^p)$ w.r.t. the
standard upper gradient to $Pf$. Together with the above inequality,
we obtain (\ref{Sobofaible}).\epr

\begin{cor}
If a Riemannian manifold $M$ with locally doubling property
satisfies $(S_{\varphi}^p)$ for the usual gradient, then it
satisfies it at any scale. If $X$ is a discretization of $M$, then
it also satisfies $(S_{\varphi}^p)$.
\end{cor}

\begin{rem}\label{finalrem}
Assume that $X$ is coarsely $b$-geodesic for every $b>0$ (e.g. $X$
is a Riemannian manifold), so that Proposition~\ref{scale} applies.
Note that in the proof of Theorem~\ref{infinitesimal->fini}, we
actually show that a Sobolev inequality at large scale is equivalent
to the Sobolev inequality for the standard upper gradient restricted
to functions of the form $g=Pf$, where $P$ is a viewpoint at some
positive scale.
\end{rem}

\section{Applications to quasi-transitive
spaces}\label{HomogeneousSection}

\subsection{Existence of a spectral gap on a quasi-transitive
metric measure space}

\subsection*{The main result}

\begin{defn}
A quasi-transitive measure space $(X,\mu)$ is a locally compact
Borel measure space on which a locally compact group $G$ acts
measurably, co-compactly, properly, and almost preserving the
measure $\mu$, i.e.
$$\sup_{g\in G}\sup_{x\in X}\frac{d(g\cdot\mu)}{d\mu}(x)<\infty.$$
\end{defn}

\begin{defn}
We call a metric measure space $(X,d,\mu)$ a quasi-transitive metric
measure space if  $(X,\mu)$ is a quasi-transitive measure space and
if $d$ is a $G$-invariant metric on $X$ which is proper and finite
on compact sets.
\end{defn}

\begin{prop}
If $G$ is $\sigma$-compact, then every $G$-quasi-transitive measure
space $(X,\mu)$ can be equipped with a metric $d$ such that
$(X,d,\mu)$ is a quasi-transitive metric measure space.
\end{prop}
\bpr Note that $d$ is not supposed to be continuous on $X$. We start
with a proper $G$-invariant metric on $G$ \cite[Theorem~7.2]{Hj}.
Take a fundamental domain $D$ in $X$ relative to the $G$-action. As
the action is co-compact, we can assume that $D$ is relatively
compact. Let $K$ be the intersection of all stabilizers of elements
of $D$. As the action is proper, $K$ is a compact subgroup of $G$.
Consider the $G$-invariant metric on $G/K$ obtained, first by
averaging our metric on $G$ over $K$ (i.e. replacing it by
$\int_{K}d(gk,hk)dk$), and then by lifting the corresponding
bi-$K$-invariant metric to $G/K$. We have a natural map $\alpha:
X\to G/K$, where $\alpha(x)$ is the unique $gK$ such that $x\in gD$.
Pulling the metric of $G/K$ to $X$ yields a $G$-invariant
pseudo-metric on $X$ which is proper and finite on compact sets. To
obtain a true metric, one can for instance add the discrete metric
on $X$ (i.e. such that two distinct points are at distance $1$).
\epr

\medskip

The following theorem is therefore more general than
Theorem~\ref{mainThmIntro}.
\begin{thm}\label{quasitransitivethm} Let $G$ be a locally compact
group and let $(X,d,\mu)$ be a quasi-$G$-transitive metric measure
space. Then $G$ is unimodular and amenable if and only if for $h$
large enough (resp. for any $h$) and every reversible viewpoint $P$
at scale $h$ on $(X,d,\mu)$, the spectral radius $\rho(P)=1,$ or in
other words, if the discrete Laplacian $\Delta=I-P$ has no spectral
gap around zero.
\end{thm}

\bpr The proof splits in three parts. First, by
Theorem~\ref{randomthm}, one checks easily that $\rho(P)=1$ if and
only if the large scale profile $j_{X,2}(t)\to\infty$ when
$t\to\infty$. Indeed, $j_{X,2}(t)\leq C$ means that $X$ satisfies a
large-scale Sobolev inequality $(S_{\varphi}^2)$ with
$\varphi(t)=C$. Thus by Theorem~\ref{randomthm}, this happens if and
only if $p_x^{2n}(x)$ has exponential decay, i.e. if and only if
$\rho(P)<1$.

Second, take a uniform left-invariant metric on $G$. The
co-compactness, properness of the $G$-action on $X$, plus the fact
that $\mu$ is almost-preserved by $G$ imply that $G$ and $X$ are
large-scale equivalent (this is straightforward). Hence, by
Theorem~\ref{largescalethm}, it is enough to prove
Theorem~\ref{quasitransitivethm} for $X=G$. This third step will be
achieved by Corollary~\ref{geomamenablecor}.

\begin{rem} Note that if we assume $G$ compactly generated, then it is
classical and not difficult to see that a quasi-$G$-transitive
metric measure space is quasi-isometric to $G$, equipped with the
word metric $d_S$ corresponding to a compact generating subset $S$
of $G$.
\end{rem}

\begin{cor}\label{riamannianSpectralGap}
Let $M$ be a Riemannian manifold and let $G$ be a closed subgroup of
$\Isom(M)$ acting co-compactly on $M$. Then $G$ is unimodular and
amenable if and only if the spectral radius of the heat kernel on
$M$ equals $1$, or in other words, if the (Riemannian) Laplacian on
$M$ has no spectral gap around zero.
\end{cor}
\bpr The Laplacian has a spectral gap if and only if $M$ satisfies a
Sobolev inequality $\|\nabla f\|_2\geq c\|f\|_2$ for the usual
gradient. As $M$ is quasi-transitive, it is easy to check that it
satisfies a local Poincaré inequality as in
Proposition~\ref{finite->infinitesimal}. Indeed, one has to prove
that such a local Poincaré inequality $(P(1,q))_{loc}$ holds, for
any $q\geq 1$ on a compact subset $K$ such that $X=\cup_{g\in G}
gK$. But this results from the fact that such inequality holds in
$\R^d$. Now, applying Proposition~\ref{finite->infinitesimal} and
Theorem~\ref{infinitesimal->fini}, we see that the spectral gap is
equivalent to a large-scale Sobolev inequality. We conclude thanks
to Theorem~\ref{quasitransitivethm}.\epr

\subsection*{Locally compact groups}

All the locally compact groups considered here are $\sigma$-compact.
Recall (see Section~\ref{ExGroups}) that a $\sigma$-compact locally
compact group can be endowed with a ``large-scale" structure of
metric measure space. Let us consider the following natural
question: is amenability a geometric property among compactly
generated locally compact groups? Recall that a locally compact
group is called amenable if it admits a left invariant mean
\cite{Pier}. By {\it geometric property}, we mean a property
characterized in terms of metric measure space. Moreover, we expect
such a property to be invariant under large-scale equivalence. F\o
lner's characterization of amenability implies that the answer is
positive when the group is {\it finitely generated}. On the
opposite, note that any connected Lie group admits a co-compact
amenable subgroup (take for instance a maximal solvable subgroup)
and therefore is always quasi-isometric to a compactly generated
locally compact amenable group. So the answer is negative in
general. Actually, we will see that the answer is yes if and only if
the group is unimodular.

Let $G$ be a $\sigma$-compact locally compact group equipped with
some proper left-invariant metric $d$ and with its Haar measure
$\mu$. Fix some $h>0$. We define the boundary of a subset $A$ of $G$
by
\begin{eqnarray*}
\partial_h A & = & A B(e,h)\cap A^cB(e,h).
\end{eqnarray*}

It is important to note that the multiplication by elements of
$B(e,h)$ is on the {\it right}, so that $AB(e,h)$ has the following
metric interpretation:
$$AB(e,h)=\cup_{x\in A} B(x,h)=[A]_h$$
where $[A]_h=\{x\in G, \ d(x,A)\leq h\}.$ In particular, this
definition of boundary coincides with the one we gave in
introduction for a general metric space.

For any sequence of compact subsets with positive measure $(F_n)$ of
$G$ and for every $g\in G$, we define
$\phi_n(g)=\mu(gF_n\vartriangle F_n)/\mu(F_n)$. Note that here, the
multiplication by $g$ is on the left.

Recall \cite{Pier} that the group $G$ is amenable if and only one of
the following equivalent statements holds:

\noindent(1) There exists a sequence $(F_n)$ such that $\phi_n(g)$
is pointwise converging to zero.

\noindent(2) There exists a sequence $(F_n)$ such that $\phi_n(g)$
converges to zero uniformly on compact sets.

\noindent(3) There exists a sequence $(F_n)$ such that $\mu(QF_n\cap
QF_n^c)/\mu(F_n)\rightarrow 0$ for every compact subset $Q$.

If a sequence $(F_n)$ satisfies (1), or equivalently, (2), then it
is called a F\o lner sequence.

\begin{rem}
Generally, in the definition of F\o lner sequence, $(F_n)$ is also
asked to be an increasing exhaustion of $G$ (this also characterizes
amenability).
\end{rem}

Here, the multiplication by $Q$ is on the {\it left}, so that
amenability is {\it not} {\it a priori} characterized in terms of
isoperimetry, or in other words, in terms of metric measured space
properties. Let us define a {\it geometric} version of amenability.

\begin{defn}
The group $G$ is called geometrically amenable if it admits a
sequence of compact subsets $(F_n)$ such that one of the following
equivalent statements holds:

\noindent(1) $\mu(F_n\vartriangle F_ng)/\mu(F_n)\rightarrow 0$ for
every $g\in G$.

\noindent(2) For every compact subset $Q$ of $G$,
$$\mu(F_nQ\cap F_n^cQ)/\mu(F_n)\rightarrow 0.$$
\end{defn}
The following proposition justifies the term ``geometric".

\begin{prop}\label{unimodulaire+moyennableprop}
A $\sigma$-compact locally compact group $G$ is geometrically
amenable if and only if for $h$ large enough, the isoperimetric
profile $j_{G,1}$ (resp. $j_{G,p}$ for any $p\geq 1$) at scale $h$
is unbounded.
\end{prop}
\bpr Clearly, (2) of the definition of geometrically amenable
implies that $j_{G,1}$ is unbounded at any scale. Conversely, the
negation of (2) together with the $\sigma$-compacity of $G$ yields
the existence of a compact subset $K$ of $G$ such that for every
measurable subset $A$ with finite measure,
$$\mu(A)\leq C\mu(AK\vartriangle A)$$
for some constant $C<\infty$. Let $h$ be such that $K\subset
B(e,h)$. It follows that
$$\mu(A)\leq C\mu(\partial_h A),$$
which means that the profile $j_{X,1}$ at scale $h$ is bounded.\epr

\

If $G$ is unimodular, up to replacing $F_n$ with $F_n^{-1}$, it is
equivalent for $G$ to have left or right F\o lner sequences. In
particular, if a group is unimodular, then it is geometrically
amenable if and only if it is amenable. Actually, we have better:
geometric amenability is equivalent to amenability plus
unimodularity.

\begin{lem}\label{nonunimodular}
If the group $G$ is non-unimodular, then it satisfies the following
isoperimetric inequality for $h$ large enough
$$\mu(\partial_h A)\geq c\mu(A)\quad \forall A\subset G$$
where $c$ is some positive constant.
\end{lem}
\bpr Let $\delta$ be the modular function of $G$. Since $G$ is
non-unimodular, there exists $g\in G$ such that $\delta(g)>1$. So,
choosing $h$ large enough, we can assume that $g\in B(e,h)$. Then
for any compact subset $A\subset G$, we have
$$\mu(\partial_h A)\geq \mu(Ag\vartriangle A)\geq \mu(Ag)-\mu(A)=(\delta(g)-1)\mu(A).\;\blacksquare$$

\begin{prop}\label{geomamenable}
Let $G$ be a $\sigma$-compact locally compact group equipped with a
left Haar measure. Then $G$ is amenable and unimodular if and only
if it admits a geometric F\o lner sequence. In particular if $G$ is
compactly generated, then $G$ is amenable and unimodular if and only
if it is geometrically amenable.
\end{prop}
\bpr This is a direct consequence of Lemma~\ref{nonunimodular} and
of the above discussion.\epr

Recall that quasi-isometries between homogeneous metric measure
spaces are large-scale equivalences. We have the following
corollaries to Theorem~\ref{largescalethm}.
\begin{cor}\label{geomamenablecor'}
Geometric amenability is invariant under large-scale equivalence
between $\sigma$-compact locally compact groups.
\end{cor}
\begin{cor}\label{geomamenablecor}
Geometric amenability is invariant under quasi-isometry between
compactly generated locally compact groups.
\end{cor}

The following corollary follows from Propositions
\ref{geomamenable}, \ref{connectivitegoupprop} and
Theorem~\ref{largescalethm}.
\begin{cor}\label{nongeomCheegerconstant}
A compactly generated locally compact group is not geometrically
amenable if and only if it is quasi-isometric to a graph with
positive Cheeger constant.
\end{cor}

\begin{cor}
Being amenable and unimodular is invariant under large-scale
equivalence between $\sigma$-compact locally compact
groups.$~\blacksquare$
\end{cor}

\bigskip
\footnotesize

\noindent \noindent Romain Tessera\\
Department of mathematics, Vanderbilt University,\\ Stevenson
Center, Nashville, TN 37240 United,\\ E-mail:
\url{tessera@clipper.ens.fr}

\end{document}